\documentclass{article}
\usepackage{amsmath,amsthm,amssymb,amsfonts}
\theoremstyle{plain}
\newtheorem{theorem}{Theorem}[section] 
\newtheorem{prop}[theorem]{Proposition}
\newtheorem{lemma}[theorem]{Lemma}
\newtheorem{corollary}[theorem]{Corollary}

\theoremstyle{remark}
\newtheorem{remark}[theorem]{Remark}

\theoremstyle{definition}

\newcommand{\zz}{{\bf Z}}
\newcommand{\qq}{{\bf Q}}

\newcommand{\cc}{{\bf C}}

\newcommand{\E}{ {\mathcal E} }      

\renewcommand{\O}{ {\mathcal O} }
\newcommand{\Q}{ {\mathcal Q} }
\newcommand{\Gal}{{\rm Gal}}        
\newcommand{\End}{ {\rm End} }      
\newcommand{\Hom}{ {\rm Hom} }      
\newcommand{\ord}{ {\rm ord} }
\newcommand{\lcm}{ {\rm lcm} }

\newcommand{\inv}{ {\rm inv} }      
\newcommand{\kernel}{{\rm Ker}}     
\newcommand{\cokernel}{{\rm Coker}} 

\def\injlim{{\varinjlim }}
\def\projlim{{\varprojlim }}

\begin{document}
\title{Group structures of elementary supersingular abelian varieties over
finite fields}
\author{Hui Zhu\\
\small \em MSRI, 1000 Centennial Drive, Berkeley, CA 94720-5070\\
{\small Email: \tt zhu@msri.org}}
\date{}
\maketitle
\begin{abstract}
Let $A$ be a supersingular abelian variety over a finite field ${\bf k}$
which is ${\bf k}$-isogenous to a power of a simple abelian variety over
${\bf k}$.  Write the characteristic polynomial of the Frobenius
endomorphism of $A$ relative to ${\bf k}$ as $f=g^e$ for a monic
irreducible polynomial $g$ and a positive integer $e$. We show that
the group of ${\bf k}$-rational points $A({\bf k})$ on $A$ is
isomorphic to $(\zz/g(1)\zz)^e$ unless $A$'s simple component is
of dimension $1$ or $2$, in which case
we prove that $A({\bf k})$ is isomorphic to
$(\zz/g(1)\zz)^a\times(\zz/{\frac{g(1)}{2}}\zz\times\zz/2\zz)^b$
for some non-negative integers $a,b$ with $a+b=e$. In particular, if
the characteristic of ${\bf k}$ is $2$ or $A$ is simple of dimension
greater than $2$, then $A({\bf k})\cong (\zz/g(1)\zz)^e$.

\end{abstract}

\section{Introduction}\label{S1}

We list some notation and terminology for this paper as follows: ${\bf
k}$ is a finite field of characteristic $p$ with $q$ elements. Let
$\bar{\bf k}$ be an algebraic closure of ${\bf k}$.  Let $A$ be an
abelian variety of dimension $d$ defined over ${\bf k}$. Let $\pi$ be
the Frobenius endomorphism of $A$ relative to ${\bf k}$ and $f$ its
characteristic polynomial.

An abelian variety over ${\bf k}$ is {\it elementary} if it is ${\bf
k}$-isogenous to a power of a simple abelian variety over ${\bf
k}$. This definition is different from that of~\cite{Waterhouse:69}
(see~\cite[page~54]{Waterhouse-Milne:71}).  An abelian variety $A$ is
elementary if and only if $f=g^e$ for some monic irreducible
polynomial $g$ over $\qq$ and some positive integer $e$. An arbitrary
abelian variety is ${\bf k}$-isogenous to a product of elementary
abelian varieties, and $f=\prod_{i=1}^{t}g_i^{e_i}$ for distinct monic
irreducible polynomials $g_i$ over $\qq$ and positive integers $e_i$.
An abelian variety $A$ over ${\bf k}$ is {\it supersingular} if each
complex root of $f$ can be written in the form $\zeta\sqrt{q}$, the
product of some root of unity $\zeta$ and the positive square root
$\sqrt{q}$.  This definition is equivalent to the standard in
literature (see section~3.2).

\begin{theorem}\label{Telementary} 
Let $A$ be an elementary supersingular abelian variety over ${\bf k}$ and
$f=g^e$ as above. Then $A({\bf k})$ is isomorphic as an abelian group to
$(\zz/g(1)\zz)^e$ except in the following cases:
\begin{enumerate}
\item[(1)] $p\equiv 3\bmod
4$, $q$ is not a square, and $A$ is ${\bf k}$-isogenous to a power of a supersingular
elliptic curve with $g=X^2+q$,
\item[(2)] $p\equiv 1\bmod 4$, $q$ is not a square,
and $A$ is ${\bf k}$-isogenous to a power of a two dimensional abelian variety
with $g=X^2-q$.
\end{enumerate}
In these two exceptional cases, there are
non-negative integers $a,b$ with $a+b=e$ such that
$$A({\bf k})\cong (\zz/g(1)\zz)^a\times
\left(\zz/{\tfrac{g(1)}{2}}\zz\times\zz/2\zz\right)^b.$$
\end{theorem}

This result is particularly striking when $p=2$ or $A$ is simple with
$d > 2$ for then $A({\bf k})\cong_{\zz}(\zz/g(1)\zz)^e$. In the latter
case $A({\bf k})$ will be either cyclic or a product of two cyclic
groups, since $e=1$ or $2$. (See Proposition~\ref{Pe12}).

We call an elementary supersingular abelian variety $A$ {\it
exceptional} if it belongs to either of the two isogeny classes stated in
Theorem~\ref{Telementary} (1) and (2). We will show
(see Proposition~\ref{Paprime}) that if $A$ is exceptional,
then for every pair of
non-negative integers $a',b'$ with $a'+b'=e$, there exists an
abelian variety $A'$ which is ${\bf k}$-isogenous to $A$
with $$A'({\bf k})\cong (\zz/g(1)\zz)^{a'}\times
\left(\zz/{\tfrac{g(1)}{2}}\zz\times\zz/2\zz\right)^{b'}.$$

In this paper $\End_{\bf k}(A)$ denotes the ring of ${\bf
k}$-endomorphisms of $A$.  Write $\End_{\bf k}^0(A)=\End_{\bf
k}(A)\otimes_{\zz}\qq$. Let $\qq[\pi]$ be the $\qq$-subalgebra of
$\End_{\bf k}^0(A)$ generated by $\pi$, let $\O$ be its maximal order,
and $\zz[\pi]$ its $\zz$-subalgebra generated by $\pi$.  The group
$A(\bar{\bf k})$ is naturally an $\End_{\bf k}(A)$-module.  Our
results describe $A(\bar{\bf k})$ as a module over any subring of
$\End_{\bf k}(A)\cap\qq[\pi]$ that contains $\zz[\pi]$.  The Galois group
$\Gal(\bar{\bf k}/{\bf k})$ is (geometrically) generated by the
Frobenius $\pi$, the $\zz[\pi]$-module structure of $A(\bar{\bf k})$
is also its Galois module structure.

For any prime number $l$ we write $R_{(l)}$ (with parenthesis) for the
localization of a commutative ring $R$ at $l$, this notation should
not be confused with $R_l$ that for the $l$-adic completion of $R$.

\begin{theorem}\label{Tmodule-elementary} 
Let $A$ be an elementary supersingular abelian variety over ${\bf k}$
of dimension $d$. Let $R$ be a ring with
$\zz[\pi]\subseteq R\subseteq \End_{\bf k}(A)\cap \qq[\pi]$.
Then there is a surjective $R$-module homomorphism
$$\varphi: A(\bar{\bf k})\rightarrow (R_{(p)}/R)^e$$
such that the cardinality of the kernel of $\varphi$ divides $2^d$.
Furthermore, $\varphi$ is an isomorphism when $p=2$.
\end{theorem}

Suppose $A$ is a simple supersingular abelian variety over ${\bf k}$
and $R$ the endomorphism ring $\End_{\bf k}(A)\cap\qq[\pi]$: if
$d\neq 2$, then $A(\bar{\bf k})\cong_R
(R_{(p)}/R)^e$; if $d=2$, then
$A(\bar{\bf k})\cong_R (R_{(p)}/R)^a\times (\O_{(p)}/\O)^b$ for some
non-negative integers $a,b$ with $a+b=e$. (See Proposition~\ref{Psimple2}.)

The group structure of the ${\bf k}$-rational points and the Galois
module structure of the $\bar{\bf k}$-rational points on an elliptic
curve were studied by~\cite{Deuring:41} (see
also~\cite[Chapter~V]{Silverman:86} and~\cite{HWL:96}). The group
structure of the ${\bf k}$-rational points on a supersingular elliptic
curve was carried out in~\cite[Chapter~4, (4.8)]{Schoof:87} (see also
Corollary~\ref{Csimple}). Our present paper yields
a description of this  nature for higher dimensional abelian varieties.
Our result for arbitrary supersingular abelian varieties
are prepared separately in~\cite{paper2}.  (Recently, independent
of our work, the group structure of dimensional two supersingular
abelian varieties was studied in~\cite{Xing:96}.)

We develop the following idea for studying the group structure of the
rational points on an elementary supersingular abelian
variety $A$ over ${\bf k}$: we show that the ring $\zz[\zeta\sqrt{q}]$
is a {\it Bass order} over some suitable subring (see
section~\ref{S5}). Next we describe the Tate modules of $A$ over
$\zz[\pi]$. Finally the group structure of $A({\bf k})$ follows by
viewing $A({\bf k})$ as the kernel of the isogeny $\pi-1$ on
$A(\bar{\bf k})$ (see section~3). Proofs of Theorems~1.1 and 1.2 lie
in section~3.

This paper is based on a portion of the author's Berkeley Ph.D thesis.
The author is deeply grateful to her advisor Professor Hendrik
Lenstra for inspiration and guidance.  The author also wish to
thank Bjorn Poonen and Phil Ryan for their comments on an earlier
version of this paper. The author was supported as a MSRI
postdoctoral fellow while preparing this paper.

\section{Torsion-free modules over Bass orders}\label{S5}

\subsection{Notions from algebra}\label{S51}

We begin this section with some notions from algebra and
then some auxiliary results from algebraic number theory.
The material largely follows~\cite[Introduction and Chapter
3]{Curtis-Reiner:90}. Here we assume all rings are commutative and
modules are {\em finitely generated}.  Let $K$ be a local or global field of
zero characteristic and $\O_K$ its discrete valuation ring or its
ring of integers, respectively.

Suppose $L$ is a finite dimensional separable $K$-algebra.  An
$\O_K$-algebra $\Lambda$ is called an {\it $\O_K$-order} (in $L$) if
it is a finitely generated projective $\O_K$-module (and
$\Lambda\otimes_{\O_K}K=L$).  A $\zz$-order is simply called an {\it
order}.  Let $\Lambda$ be an $\O_K$-order in $L$. We denote the unique
maximal $\O_K$-order in $L$ by $\O_L$.  If $M$ is a $\Lambda$-module
which is projective over $\O_K$, then $M$ is called a {\it
$\Lambda$-lattice}.

For any prime $\wp$ of $\O_K$, we denote by $(\O_K)_\wp, \Lambda_\wp,
M_\wp$ their $\wp$-adic completions, respectively. If $K=\qq$ and
$\wp = l$ for some prime number $l$, then we write $(\O_K)_l, \Lambda_l,
M_l$ for their $l$-adic completions.

A $\Lambda$-module $M$ is called {\it torsion-free} if $\alpha m\neq 0$ for
any non-zerodivisor $\alpha\in \Lambda-\{0\}$ and $m\in
M-\{0\}$. In particular, when $\Lambda$ is a domain
then this is equivalent to the standard definition.
A $\O_K$-module is projective if and only if it is
torsion-free~\cite[II.4 (4.1)]{Frohlich-Taylor:91}.
So $M$ is a torsion-free $\Lambda$-module
if and only if it is a $\Lambda$-lattice.  If $M$ is a torsion-free
$\Lambda$-module, then it is torsion-free over $\O_K$, hence there is
a natural embedding $M\hookrightarrow M\otimes_{\O_K}K$, where
$M\otimes_{\O_K}K$ has a natural $L$-module structure.
If $M\otimes_{\O_K}K$ is free of rank $e$ over $L$ for some integer $e$,
then $M$ is said of {\it rank} $e$.
We shall note here that
$e$ is used to denote an arbitrary positive integer in this section.

Suppose $L$ is a finite field extension of $K$.
Denote by $\Delta_{\Lambda/\O_K}$ the
discriminant (ideal) of $\Lambda$ over $\O_K$
and $\Delta_{L/K}:=\Delta_{\O_L/\O_K}$.
We recall that
$[\O_L:\Lambda]^2\,\Delta_{L/K}=\Delta_{\Lambda/\O_K}$ and so
$[\O_L:\Lambda]^2\mid\Delta_{\Lambda/\O_K}.$
Let $\alpha$ be an integral element in $L$ and $h\in \O_K[X]$ be its (monic) minimal
polynomial.
Let $\Lambda=\O_K[\alpha]$. Then
$\Delta_{\Lambda/\O_K}=\O_K\,\Delta(h)$.  Let $\wp$ be any non-zero
prime ideal of $\O_K$.  Then $\Lambda_\wp$ is semilocal and
$\Lambda_\wp\cong \prod_{Q\mid \wp}\Lambda_Q$ where the product
ranges over all prime ideals $Q$ of $\Lambda$ lying over $\wp$.  There
is a bijective correspondence between these $Q$'s of $\Lambda$ and the
set of (monic) irreducible factors $\bar{h_0}$ of $\bar{h}=(h\bmod
\wp)\in ({\O_K}/\wp)[X]$.  (See~\cite[Chapter I, Proposition~25, page
27]{Lang:86}.)  If $Q$ corresponds to $\bar{h_0}$ in this bijection,
then $Q=(\wp,h_0(\alpha))$ in $\Lambda$.  We use the following
notation throughout this paper: for any prime ideal $v$ of $\O_L$
lying over a prime $\wp$ of $\O_K$, let $\gamma(v/\wp)$,
$\kappa(v/\wp)$ and $\varrho(v/\wp)$ denote the {\em ramification
index}, {\em residue field degree} and {\em decomposition degree},
respectively. In particular, when $\Lambda = \O_L$ then
$\kappa(Q/\wp)=\dim_{{\O_K}/\wp}\Lambda/Q=\deg(\bar{h_0})$ and
$\gamma(Q/\wp)$ equals the multiplicity of $\bar{h_0}$ as a factor of
$\bar{h}$.  We have the following fundamental lemma.  (This proof is
due to Hendrik Lenstra.)

\begin{lemma}\label{Llenstra}
Let the notation be as above. Then the prime ideal $Q$ is not invertible if
and only if $\bar{h_0}^2\mid\bar{h}$ and all coefficients of the
remainder of $h$ upon division by $h_0$ are in $\wp^2$.  The
$(\O_K)_\wp$-order $\Lambda_\wp$ is not the maximal order if and only
if there is a monic irreducible factor $\bar{h_0}$ of $\bar{h}$ such
that $\bar{h_0}^2\mid\bar{h}$, and all coefficients of the remainder
of $h$ upon division by $h_0$ are in $\wp^2$.
\end{lemma}
\begin{proof}
Write $J: = (\wp,h_0(X))$ in $\O_K[X]$, it is a prime ideal.  The
natural surjective map $\O_K[X]\rightarrow \Lambda$ induces a
surjective map $\theta: J/J^2\rightarrow Q/Q^2$ with $\kernel(\theta)$
generated by $h$. Write $h$ in base $h_0$ and obtain
$h=r_2h_0^2+r_1h_0+r_0$ for some $r_2,r_1,r_0\in {\O_K}[X]$ with
$\deg(r_1), \deg(r_0)<\deg(\bar{h_0})$. Then $h\in J$ if and only if
$r_0\in \wp[X]$, while $h\in J^2$ if and only if $r_1\in \wp[X]$ and
$r_0\in \wp^2[X]$. So
$\dim_{\Lambda/Q}J/J^2=1+\dim_{{\O_K}/\wp}\wp/\wp^2=2$ and hence
$\dim_{\Lambda/Q}Q/Q^2=\dim_{\Lambda/Q}J/J^2-\dim_{\Lambda/Q}\kernel(\theta)
=2-\dim_{\Lambda/Q}\kernel(\theta)$, where
$\dim_{\Lambda/Q}\kernel(\theta)$ is $0$ or $1$.  Therefore,
$\dim_{\Lambda/Q}Q/Q^2\neq 1$ if and only if $h\in J^2$. We conclude
that $Q$ is not invertible if and only if $h\in J^2$, which is
equivalent to $\bar{h_0}^2\mid \bar{h}$ and all coefficients of the
remainder of $h$ upon division by $h_0$ are in $\wp^2$.  Thus the
semilocal ring $\Lambda_\wp$ is maximal if and only if $\Lambda_Q$ is
maximal for each prime ideal $Q$ over $\wp$, which is equivalent to
$Q$ is invertible, and so follows our assertion.
\end{proof}

\begin{corollary}\label{Clenstra}
Let the notation be as in Lemma~\ref{Llenstra}.
If $h_0=X-\beta$ with $\beta\in \O_K$, then
$Q$ is not invertible if and only if $h(\beta)\equiv 0\bmod \wp^2$
and $h'(\beta)\equiv 0\bmod \wp$, where $h'$ denotes
the derivative of $h$.
\end{corollary}
\begin{proof}
The condition $\bar{h_0}^2\mid \bar{h}$ is equivalent to that
$h(\beta)\equiv 0\bmod \wp$ and
$h'(\beta)\equiv 0\bmod \wp$.
The condition that all coefficients of the remainder of
$h$ upon division by $h_0$ are in
$\wp^2$ is equivalent to $h(\beta)\equiv 0\bmod \wp^2$.
\end{proof}

\subsection{Bass orders}

A reference for concepts in this subsection
is~\cite[Chapter~4]{Curtis-Reiner:90}.  Let $K$ and $\O_K$ be as the
previous subsection.  Let $L$ be a finite field extension over $K$.
We call an $\O_K$-order $\Lambda$ a {\it Gorenstein order} if every
exact sequence of $\Lambda$-modules $0\rightarrow \Lambda\rightarrow
M\rightarrow N\rightarrow 0,$ in which $M$ and $N$ are
$\Lambda$-lattices is split over $\Lambda$. If $\Lambda$ has the
additional property that every $\O_K$-order in $L$ containing
$\Lambda$ is also a Gorenstein order, then we call $\Lambda$ a {\it
Bass order}. Note that being a Bass order is a local property, in
other words, $\Lambda$ is a Bass $\O_K$-order if and only if
$\Lambda_\wp$ is a Bass $(\O_K)_\wp$-order for every prime $\wp$ in
$\O_K$.

\begin{prop}\label{Pequivalent}
The following are equivalent:
\begin{enumerate}
\item[(1)] $\Lambda$ is a Bass $\O_K$-order;
\item[(2)] $\O_L/\Lambda$ is a cyclic $\Lambda$-module;
\item[(3)] every ideal of $\Lambda$ can be generated by two elements;
\item[(4)] for every maximal ideal $Q$ of $\Lambda$ we have
$\dim_{\Lambda_Q/Q\Lambda_Q}(\O_L)_Q/Q(\O_L)_Q \leq 2$;
\item[(5)] the multiplicity of $\Lambda$ at each maximal ideal $Q$ is $\leq 2$.
\end{enumerate}
\end{prop}
\begin{proof}
The first three parts are equivalent according
to~\cite[Theorem~2.1]{Levy:85}.
The last two parts are equivalent to (1) by~\cite[Theorem~2.1]{Greither:82}.
\end{proof}

\begin{remark}\label{Xbass}
Here are some examples of Bass orders of interest.
\begin{enumerate}
\item[(i)]
If $L$ is a quadratic field extension over $K$, then
$(\O_L)_\wp/\Lambda_\wp$ is cyclic over $\Lambda_\wp$ for
every prime $\wp$ of $\O_K$ and thus $\Lambda$ is a Bass order over $\O_K$.
\item[(ii)] All maximal orders in number fields are Bass orders.
\end{enumerate}
\end{remark}

We are interested in describing torsion-free modules $M$ over $\Lambda_\wp$
of rank $e$ for a prime ideal $\wp$ of $\O_K$.
Recall that $\Lambda_\wp$ is a semilocal ring whose
maximal ideals are those prime ideals $Q$ lying over $\wp$, so there is
a corresponding decomposition of $M$ as $M\cong \prod_{Q\mid \wp} M_Q$.
If $\Lambda_\wp$ is maximal, that is
$\Lambda_\wp =(\O_L)_\wp$, then $M_Q$ is
torsion-free over the principal ideal domain $(\O_L)_Q$ of rank $e$,
so $M_Q\cong\Lambda_Q^e$ for all $Q$. Thus $M\cong\Lambda_\wp^e$.
If $\Lambda_\wp$ is not maximal, then it is generally hard to
classify such modules $M$ in terms of orders in $L_\wp$
(see~\cite[Chapter~3]{Curtis-Reiner:90}).
However, torsion-free modules over Bass orders can be described as follows.

\begin{theorem}[Bass]\label{Tbass}
Let $K$ be a local field, $\O_K$ its discrete valuation ring,
and $\Lambda$ a Bass $\O_K$-order in a finite field extension $L$ over $K$.
Then every indecomposable
torsion-free $\Lambda$-module is a projective $\Lambda'$-module for some
$\O_K$-order $\Lambda'$ in $L$ containing $\Lambda$.
\end{theorem}
\begin{proof}
Follows from the equivalencies in Proposition~\ref{Pequivalent},
~\cite[Theorem (37.13)]{Curtis-Reiner:90}
and the definition of Bass orders.
\end{proof}

\subsection{Supersingular $q$-numbers}\label{S3}

This subsection contains a technical part of this paper, which lies in
Lemma~\ref{Lindex}. We first of all introduce some notations.
For any positive integer $n$, and any prime number $l$,
let $n_l$ and $n_{(l)}$ denote the $l$-part and
the non-$l$-part of $n$ respectively;
let $\zeta_n=\exp(\frac{2{\rm\pi}\sqrt{-1}}{n})$.
The primitive $n$-th roots of unity are the $\zeta_n^\nu$
with positive integers $\nu$ that are coprime to $n$.

For the rest of the paper $l$ is a prime number different from $p$.
An algebraic number $\alpha\in \cc$ is called a {\it supersingular
$q$-number} if it is of the form $\zeta\sqrt{q}$.  (See section~3.2
for its relationship to supersingular abelian variety.)  Write
$\pi=\zeta_m^{\nu}\sqrt{q}$.  Let $K=\qq(\pi^2)$ and let $\O, \O_K$ be
the ring of integers of $\qq(\pi),K$, respectively.  Obviously
$K=\qq(\zeta_{m/(2,m)})$ and $[\qq(\pi):K]=1$ or $2$.  In this paper,
we write $(n_1,n_2)$ for the greatest common divisor of two integers
$n_1, n_2$, we denote by $(\frac{\cdot}{p})$ the Jacobi symbol.

To prove the following two lemmas we need a few well-known and
elementary results from algebraic number theory, which  we
recall here for the convenience
of the reader: For any prime number $p$ and
positive integer $n$ we have (1) $\Delta_{\qq(\sqrt{p})/\qq}=p$ if
$p\equiv 1\bmod 4$, and $4p$ if $p\not\equiv 1\bmod 4$; (2)
$\sqrt{p}\in \qq(\zeta_n)$ if and only if
$\Delta_{\qq(\sqrt{p})/\qq}\mid n$; (3) Let $p\neq 2$, if $p\mid n$
then $\qq(\sqrt{(\frac{-1}{p})p})\subseteq\qq(\zeta_n)$.

\begin{lemma}\label{LRF}
Suppose $q$ is a non-square. Then  $\qq(\pi)=K$ if and only if
\begin{eqnarray*}
(1)\ \Delta_{\qq(\sqrt{p})/\qq}\mid m,\quad\mbox{and}\quad
(2)\ \Delta_{\qq(\sqrt{p})/\qq}\nmid m/(2,m)\mbox{ if }4\mid m.
\end{eqnarray*}
In this case $2\mid \gamma(v/p)$ for any prime $v$ of $\O$ lying over $p$.
\end{lemma}
\begin{proof}
We note $\qq(\pi)=\qq\left(\zeta_{m/(2,m)},\sqrt{p\zeta_{m/(2,m)}}\right)
=K\left(\sqrt{p\zeta_{m/(2,m)}}\right).$
Thus $\qq(\pi)=K$ if and only if
$\sqrt{p\zeta_{m/(2,m)}}\in K$, and if and only if
$K\left(\sqrt{p}\right)=K\left(\sqrt{\zeta_{m/(2,m)}}\right),$
that is,
$\qq\left(\zeta_{m/(2,m)},\sqrt{p}\right)=\qq(\zeta_m)$. This is equivalent to
$$(1a)\ \sqrt{p}\in\qq(\zeta_m),\quad\mbox{and}\quad (2a)\
\sqrt{p}\notin\qq(\zeta_{m/(2,m)}) \mbox{ if } 4\mid m,$$ which is
equivalent to (1) and (2) respectively by the remark preceding this lemma.

To show the second assertion it is enough to prove it for just one
prime $v$ over $p$ since all primes lying over $p$ are conjugate as
$\qq(\pi)$ is the cyclotomic field $\qq(\zeta_{m/(2,m)})$.  We claim
$(2,p)p\mid\frac{m}{(2,m)}$. In fact, if $p=2$ then (1) implies $8\mid
m$ by the remark preceding this Lemma, so our claim follows; if $p\neq
2$ then (1) implies $p \mid m$. But since $p\neq 2$, we have $p\mid
\frac{m}{(2,m)}$. By the remark preceding this lemma, we thus see that
$\qq(\zeta_{m/(2,m)})$ contains a quadratic field
$\qq(\sqrt{(\frac{-1}{p})p})$ over $\qq$ where $p$ is totally
ramified. Hence $2\mid \gamma(v/p)$.
\end{proof}

Let $\E$ be the set of supersingular $q$-numbers
$\zeta_m^{\nu}\sqrt{q}$ which satisfy the following conditions:
$p\neq 2$, $q$ is not a square, $p\nmid m$, and
\begin{eqnarray*}
\mbox{(1) $4\nmid m$ when $p\equiv
1\bmod 4$;  and  (2) $4||m$ when $p\equiv 3\bmod 4$.}
\end{eqnarray*}

For the proof of the lemma below, we remark here
that $\alpha\in \Lambda_\wp$
is a unit if and only if $\alpha$ is coprime to $\wp$.

\begin{lemma}\label{Lindex}
Let the notation be as above.
If $(l,\pi)\notin \{2\}\times\E$ then
$\zz[\pi]_l=\O_l$. If $(l,\pi)\in
\{2\}\times\E$ then $\zz[\pi]_2\subsetneq \O_2$;
let $\wp$ be any prime ideal in $\O_K$ lying over $2$, then
\begin{enumerate}
\item[(1)] $\qq(\pi)$ is a quadratic extension over $K$ where $\wp$ is split
if $p\equiv \pm 1\bmod 8$, and $\wp$ is inert if $p\equiv \pm 3\bmod 8$.
\item[(2)] \begin{sloppypar}
$\zz[\pi]_\wp$ is a local ring and a Bass $(\O_K)_\wp$-order such that
$\O_\wp$ is the only $(\O_K)_\wp$-order in $\qq(\pi)_\wp$ that
properly contains $\zz[\pi]_\wp$. Moreover,
$\O_\wp/\zz[\pi]_\wp\cong_{(\O_K)_\wp}(\O_K)_\wp/\wp$.\end{sloppypar}
\end{enumerate}
\end{lemma}
\begin{proof}
If $q$ is a square, then $\pi\notin \E$ and
$\zz[\pi]_l=\zz[\zeta_m]_l=\O_l$. For the rest of the proof, we
assume $q$ is not a square. We consider the following two cases.

\noindent{\fbox{Case 1. $l\neq 2$ or $p\mid m$.}}

We claim $\zz[\pi]_l=\O_l$.  Since $l\neq p$, we note that
$\zz[\pi^2]_l=\zz[\zeta_{m/(2,m)}]_l=(\O_K)_l$.  Suppose $l\neq 2$,
obviously $\zz[\pi^2]_l=(\O_K)_l\subseteq \zz[\pi]_l\subseteq \O_l$.
If $\qq(\pi)=K$ then $\O_l=(\O_K)_l$ and so $\zz[\pi]_l=\O_l$.  If
$\qq(\pi)\neq K$ then $[\O:\zz[\pi]]_l^2\mid
\Delta_{\zz[\pi]/(\O_K)};$ but since $\zz[\pi]\cong
\O_K[X]/(X^2-q\zeta_{m/(2,m)})$, we have
$\Delta_{\zz[\pi]/\O_K}=\O_K\Delta(X^2-q\zeta^\nu_{m/(2,m)}) =4q\O_K.$
So $(\Delta_{\zz[\pi]/\O_K})_l=(\O_K)_l$ since $4q$ is coprime to $l$.
Therefore $[\O:\zz[\pi]]_l$ is the unit ideal and so
$\zz[\pi]_l=\O_l$.

Now let $l=2$ and $p\mid m$.
By the remark preceding Lemma~\ref{LRF},
$\sqrt{(\frac{-1}{p})p}\in\zz[\zeta_{m/(2,m)}]_2=\zz[\pi^2]_2
\subseteq \zz[\pi]_2$.
Moreover, the norm of $\sqrt{(\frac{-1}{p})p}$ over $\qq$
is $\pm p$ which is coprime to $2$ so $\sqrt{(\frac{-1}{p})p}$
is a unit in $\zz[\pi]_2$.
Therefore,
$\zz[\pi]_2=\zz[\pi\sqrt{(\frac{-1}{p})p}]_2
=\zz[\zeta^\nu_m\sqrt{(\frac{-1}{p})}]_2$.  This proves our claim.

\noindent {\fbox{Case 2. $l=2$ and $p\nmid m$.}}

Write $m=2^jm_{(2)}$.  It is easy to verify that
$\qq(\pi)=\qq(\zeta_{m_{(2)}},\alpha)$ where
$\alpha=\zeta^\mu_{2^j}\sqrt{p}$ for some $2^j$-th primitive root of
unity $\zeta^\mu_{2^j}$.  We note that $\qq(\zeta_{m_{(2)}})$ and
$\qq(\alpha)$ are linearly disjoint and that the minimal polynomial of
$\alpha$ over $\qq(\zeta_{m_{(2)}})$ is $h=X^{2^{j-1}}+p^{2^{j-2}}$ if
$j\geq 2$, and is $h=X^2-p$ if $j<2$.

Let $\wp'$ be any prime ideal in the ring of integers of
$\qq(\zeta_{m_{(2)}})$ lying over $2$.
We show $\zz[\pi]_{\wp'}=\zz[\zeta_{m_{(2)}},\alpha]_{\wp'}$.
The inclusion $\zz[\pi]_{\wp'}\subseteq \zz[\zeta_{m_{(2)}},
\alpha]_{\wp'}$ is trivial.
Conversely, since $\pi^{2^j}=\zeta_{m_{(2)}}^{2^j}q^{2^{j-1}}$ and
$\alpha=\pi\zeta_{m_{(2)}}^{-\mu}$, we have $\zeta_{m_{(2)}}, \alpha\in
\zz[\pi]_{\wp'}$. Thus $\zz[\zeta_{m_{(2)}},\alpha]_{\wp'}\subseteq
\zz[\pi]_{\wp'}$.
That is, $\zz[\pi]_{\wp'}=\zz[\zeta_{m_{(2)}},\alpha]_{\wp'}$.
Hence, $\zz[\pi]_{\wp'}=(\zz[\zeta_{m_{(2)}}]_{\wp'})[\alpha]$.

If $j\geq 2$, then $h\equiv (X-1)^{2^{j-1}} \bmod \wp'$.  Note that
$\zz[\zeta_{m_{(2)}}]_{\wp'}$ is a complete
discrete valuation ring, so we have
by Corollary~\ref{Clenstra} that $\zz[\pi]_{\wp'}$ is not maximal if
and only if $h(1)=1+q^{2^{j-2}}\equiv 0\bmod {\wp'}^2$, that is,
$j=2$ and $p\equiv 3\bmod 4$. Similarly, if $j<2$ then $h\equiv
(X-1)^2\bmod \wp'$ and so $\zz[\pi]_{\wp'}$ is not maximal if and only
if $p\equiv 1\bmod 4$.  Note $\zz[\pi]_2=\prod_{\wp'\mid
2}\zz[\pi]_{\wp'}$.  By Lemma~\ref{Llenstra} and the above argument,
$\zz[\pi]_2$ is not maximal if and only if $\pi\in\E$.

In the special case $\pi\in\E$, we have
$K=\qq(\zeta_{m_{(2)}})$ and
$\qq(\pi)=K(\sqrt{(\frac{-1}{p})p})$ is quadratic over $K$.
Moreover, $\zz[\pi]_\wp = (\O_K)_\wp[\sqrt{(\frac{-1}{p})p}]$ and
$\wp$ is totally ramified in $\zz[\pi]_\wp$. This proves
that $\zz[\pi]_\wp$ is a local ring.
The decomposition of $\wp$ in the quadratic extension $\qq(\pi)$ over $K$
corresponds to that of $2$ in $\qq(\sqrt{(\frac{-1}{p})p})$ over
$\qq$, which is as in our assertion.  Since $\zz[\pi]_\wp$ is a
quadratic order over the complete discrete valuation ring
$(\O_K)_\wp$, it is a Bass order by Remark~\ref{Xbass}~(1).  As
$(\O_K)_\wp$-orders, $\zz[\pi]_\wp
\subset\O_\wp\cong (\O_K)^2_\wp$.
There is an injection
$\zz[\pi]_\wp/(\O_K)_\wp\hookrightarrow \O_\wp/(\O_K)_\wp\cong
(\O_K)_\wp$, under which $\zz[\pi]_\wp/(\O_K)_\wp\cong
\wp^i (\O_K)_\wp$ for some positive integer $i$.
In other words, $\zz[\pi]_\wp=(\O_K)_\wp+\wp^i\O_\wp$ and so
$\O_\wp/\zz[\pi]_\wp\cong (\O_K)_\wp/\wp^i$.  But
$\Delta_{\zz[\pi]_\wp/(\O_K)_\wp}
=(\O_K)_\wp\Delta(X^2-(\frac{-1}{p})p)=4(\O_K)_\wp$ and hence
$[\O_\wp:\zz[\pi]_\wp]^2=[(\O_K)_\wp:\wp^i]^2=2^{2i}(\O_K)_\wp\mid
4(\O_K)_\wp$. Thus $i=1$, that is,
$\O_\wp/\zz[\pi]_\wp\cong(\O_K)_\wp/\wp$ as $(\O_K)_\wp$-modules.
Hence $\O_\wp$ is
the only $(\O_K)_\wp$-order in $\qq(\pi)_\wp$ that properly contains
$\zz[\pi]_\wp$.
\end{proof}

\subsection{Torsion-free modules over Bass orders}\label{S54}

Let the notation be as in section~\ref{S3}. For any ring $R$ we use
$R^*$ to denote its group of units.
Henceforth in this section we
assume that $R$ is an order in $\qq(\pi)$ containing $\zz[\pi]$.  Let $M$
be a torsion-free $R_l$-modules
(as defined in section~\ref{S51}) of
rank $e$, our goal here is to describe
all such modules. We recall that all modules are
assumed to be finitely generated.

\begin{lemma}\label{Lindex2}
Let $\wp$ be any prime ideal in $\O_K$ lying over $2$.  Let $N$ be an
indecomposable torsion-free $\zz[\pi]_\wp$-module. Suppose $(l,\pi)\in
\{2\}\times\E$.  If $\wp$ is inert in $\qq(\pi)$, then $N\cong
\zz[\pi]_\wp$ or $\O_\wp$. If $\wp$ is split, i.e., $\wp=\wp_1\wp_2$
for some prime ideals $\wp_1$, $\wp_2$ in $\qq(\pi)$, then $N\cong
\zz[\pi]_\wp$, $\O_{\wp_1}$, or $\O_{\wp_2}$.
\end{lemma}
\begin{proof}
By Lemma~\ref{Lindex}, we know that $\zz[\pi]_\wp$ is a local ring and
an $(\O_K)_\wp$-order, so we invoke Theorem~\ref{Tbass}.  If $N$ is
projective over the local ring $\zz[\pi]_\wp$ then $N\cong
\zz[\pi]_\wp$. Otherwise, $N$ is projective over $\O_\wp$, since
$\O_\wp$ is the only $(\O_K)_\wp$-order of $\qq(\pi)_\wp$ that
properly contains $\zz[\pi]_\wp$ by Lemma~\ref{Lindex}~(2).  Suppose
$\wp$ is inert in $\qq(\pi)$, that is, $\O_\wp$ is a discrete
valuation ring then $N\cong_{\O_\wp} \O_\wp$. If $\wp$ splits into
$\wp_1$ and $\wp_2$ in $\qq(\pi)$, that is, if
$\O_\wp\cong\O_{\wp_1}\times \O_{\wp_2}$, then $N\cong_{\O_\wp}
\O_{\wp_1}$ or $\O_{\wp_2}$.  Therefore
$$N\cong_{\zz[\pi]_\wp}\zz[\pi]_\wp,
\O_{\wp_1},\mbox{  or}\quad\O_{\wp_2}.$$This finishes the proof.
\end{proof}

\begin{prop}\label{Pindecomposable}
There is the following isomorphism of $R_l$-modules:
$$M\cong_{R_l}
\left\{
\begin{array}{ll}
R_l^e&\mbox{if $(l,\pi)\notin \{2\}\times\E$},\\
\prod_{\wp\mid l} (R_\wp^{a_\wp} \times
\O_\wp^{b_\wp})&\mbox{if $(l,\pi)\in \{2\}\times\E$}
\end{array}\right.$$
where $\wp$ ranges over all prime ideals in $\O_K$ lying over $2$, and
$a_\wp$, $b_\wp$ are non-negative integers such that $a_\wp+b_\wp=e$.
\end{prop}
\begin{proof}
Suppose $(l,\pi)\notin \{2\}\times\E$.  By Lemma~\ref{Lindex}, the
$\zz_l$-order $R_l$ is maximal and our assertion follows from
the argument preceding Theorem~\ref{Tbass}.

Suppose $(l,\pi)\in \{2\}\times\E$.  Since $M_\wp$ is a torsion-free
$R_\wp$-module of rank $e$, by the Krull-Schmidt-Azumaya
theorem~\cite[Theorem (30.6)]{Curtis-Reiner:90}, $M_\wp$ can be
expressed as a finite direct sum of indecomposables with the summands
unique up to isomorphism and order of occurrence.  If $\wp$ is inert
in $\qq(\pi)$, then by Lemma \ref{Lindex2} there are non-negative
integers $a_\wp, b_\wp$ with $a_\wp+b_\wp=e$ such that $M_\wp\cong
R_\wp^{a_\wp}\times \O_\wp^{b_\wp}$. Now suppose $\wp$ is split in
$\qq(\pi)$. Then $M_\wp\cong R_\wp^{a_{\wp}}\times
\O_{\wp_1}^{b_\wp}\times \O_{\wp_2}^{c_\wp}$ for some non-negative
integers $a_\wp, b_\wp, c_\wp$; by comparing ranks in
$\qq(\pi)_\wp^e\cong_{\qq(\pi)_\wp} \qq(\pi)_\wp^{a_\wp}\times
\qq(\pi)_{\wp_1}^{b_\wp}\times \qq(\pi)_{\wp_2}^{c_\wp}$, we are
forced to have $b_\wp=c_\wp$.  Thus, $M_\wp\cong R^{a_\wp}_\wp \times
(\O_{\wp_1}\times \O_{\wp_2})^{b_\wp}\cong R_\wp^{a_\wp}\times
\O_\wp^{b_\wp}$ for $a_\wp, b_\wp$ with $a_\wp+b_\wp=e$. Therefore
$$M\cong\prod_{\wp\mid 2} M_\wp \cong_{R_2}\prod_{\wp\mid 2}
(R_\wp^{a_\wp} \times \O_\wp^{b_\wp}).$$ This finishes our
proof. \end{proof}

The following corollary is prepared for the next section.

\begin{corollary}\label{C58}
If $M$ is a torsion-free $R_l$-module of rank $e$ then we have
$M/(\pi-1)M\cong_{R_l} (R_l/(\pi-1))^e$
unless $l=2$, $q$ is not a square, and
$\pi=\pm\sqrt{(\frac{-1}{p})q}$, in which case
there are non-negative integers $a,b$ with $a+b=e$ such that
$$M/(\pi-1)M \cong_{R_2} (R_2/(\pi-1))^a \times (\O_2/(\pi-1))^b.$$
\end{corollary}
\begin{proof}
First of all we show that $m\notin 2^\zz$ if and only if
$R_2/(\pi-1)=0$, that is, $\pi-1\in R_2^*$. Indeed, $m\notin 2^\zz$
implies $\zeta_m^\nu-1\in \zz[\pi]_2^*\subseteq R_2^*$.  Write
$\pi-1=(\zeta_m^\nu-1)\sqrt{q} + (\sqrt{q}-1)$.  If $p=2$ then
$(\zeta_m^\nu -1)\sqrt{q}$ lies in a prime over $2$ while
$\sqrt{q}-1\in R_2^*$ so their sum lies in $R_2^*$; if $p\neq 2$, then
$R_2^*\sqrt{q}=R_2^*$ and $\sqrt{q}-1$ lies in a prime over $2$ thus
their sum also lies in $R_2^*$. This proves our claim.  Consequently,
if $m\in 2^\zz$ then $M/(\pi-1)M\cong (R_2/(\pi-1))^e$ since they are
both trivial.  By Proposition~\ref{Pindecomposable}, we have
$M/(\pi-1)M\cong_{R_l}(R_l/(\pi-1))^e$ unless $l=2, \pi\in\E$ and
$m\in 2^\zz$.  By the definition of $\E$, we have
$\pi=\zeta_m^\nu\sqrt{q}\in\E$ if and only if $q$ is not a square and
$m=1$ or $2$ if $p\equiv 1\bmod 4$; while $m=4$ if $p\equiv 3\bmod 4$.
That is, we have $l=2$, $q$ is not a square and
$\pi=\pm\sqrt{(\frac{-1}{p})q}$.
\end{proof}

\section{Supersingular abelian varieties}\label{S7}

\subsection{Preliminaries}\label{S2}

This subsection provides some auxiliary results on abelian varieties over
finite fields. We shall quote from~\cite{Milne}
and~\cite{Mumford:74} without comment.

Recall that $l$ is any prime different from $p$.  If $G$ is an abelian
group we denote by $G[l^\infty]$ the subgroup of all elements in $G$
whose order is a $l$-power. For every ${\bf k}$-isogeny
$r:A\rightarrow A$, we denote by $A[r]$ the kernel of the induced map
on $A(\bar{\bf k})$ as abelian groups.  The $l$-adic Tate module
$T_l(A)=\projlim_n A[l^n]$ is free of rank $2d$ over $\zz_l$. Since
the Frobenius endomorphism $\pi$ acts faithfully on it, $T_l(A)$ is
torsion-free $\zz[\pi]_l$-module, and $V_l(A):=T_l(A)\otimes_{\zz_l}\qq_l$
is a $\qq[\pi]_l$-module.  We also know that $\qq[\pi]$ is a
semisimple $\qq$-algebra. If the characteristic polynomial of the
Frobenius is $f=\prod_{i=1}^{t}g_i^{e_i}$ as in section~1, then
\begin{eqnarray*}
\qq[\pi]_l  \cong \prod_{i=1}^{t} \qq[\pi]_l/(g_i(\pi)),\quad
V_l(A) \cong \prod_{i=1}^{t}
(\qq[\pi]_l/(g_i(\pi)))^{e_i}.
\end{eqnarray*}
In particular,
if $A$ is elementary so $\qq[\pi]\cong\qq[\pi]/(g(\pi))$ is a field,
and we note that $V_l(A)\cong\qq(\pi)_l^e$.
Thus $T_l(A)$ is a torsion-free module of rank $e$ over
any $\zz_l$-order of $\qq(\pi)_l$ containing $\zz[\pi]_l$.

It is known that $T_l$ defines a (covariant) functor from the category
of abelian varieties $A'$ over ${\bf k}$ with a ${\bf k}$-isogeny $r:
A\rightarrow A'$ to the category of $\zz[\pi]_l$-lattices (as
$\zz_l$-order) $T_l(A')$ of $V_l(A)$ with an injective
$\zz[\pi]_l$-module homomorphism $r: T_l(A)\rightarrow T_l(A')$. In
fact, every $\zz[\pi]_l$-lattice of $V_l(A)$ containing $T_l(A)$
arises this way (see Proposition~\ref{P21}).  Note $V_l(A)/T_l(A)\cong
A[l^\infty]$. Mapping the short exact sequence $0\rightarrow
T_l(A)\rightarrow V_l(A)\rightarrow A[l^\infty]\rightarrow 0$ to that
of $A'$ by $r$ induces an injective $\zz[\pi]_l$-module homomorphism
$r: T_l(A)\rightarrow T_l(A')$ with cokernel $T_l(A')/r T_l(A)$ and an
isomorphism $V_l(A)\rightarrow V_l(A')$. Let $r^{-1}T_l(A')$ be the
pullback of $T_l(A')\subset V_l(A')$ under this isomorphism, there is
an isomorphism $T_l(A')/r T_l(A)\cong r^{-1}T_l(A')/T_l(A)$. Applying
the Snake Lemma to the above resulting diagram, we have
$r^{-1}T_l(A')/T_l(A)\cong\kernel(r)[l^\infty]$,
where $\kernel(r)$
denotes the kernel (as abelian groups) of the induced map
$A(\bar{\bf k})\stackrel{r}{\rightarrow}
A'(\bar{\bf k})$.

\begin{prop}\label{P21}
For any prime $l\neq p$, let $\theta:
V_l(A)/T_l(A)\stackrel{\sim}{\rightarrow} A[l^\infty]$ be the
isomorphism as above.  For every $\zz[\pi]_l$-lattice $M$ containing
$T_l(A)$, there is an abelian variety $A'$ with a ${\bf k}$-isogeny
$r: A\rightarrow A'$ such that $M=r^{-1}T_l(A')$ in $V_l(A)$ and
$\theta(M/T_l(A))=\kernel(r)$.
\end{prop}
\begin{proof}
Write $G:=\theta(M/T_l(A))$.  We note that $G$ is a finite subgroup of
$A(\bar{\bf k})$ of $l$-power order (coprime to $p$) and it has a
induced $\Gal(\bar{\bf k}/{\bf k})$-module structure.  So it
determines a finite \'etale subgroup scheme ${\cal G}$ of $A$ over ${\bf k}$
with ${\cal G}(\bar{\bf k})=G$.  Take $A'= A/{\cal G}$ and the obvious
${\bf k}$-isogeny $r:A\rightarrow A'$, we see that $\kernel(r)=G$.
The argument preceding the proposition indicates that
$\theta(r^{-1}T_l(A')/T_l(A))=G$. Our assertion follows.
\end{proof}

Define $T_p(A)=\projlim_nA[p^n]$ in an analogous manner.
It is free $\zz_p$-module of rank between $0$ and $d$ (inclusive).
(There is more on this in section~3.2.)

To begin our study of the group structure of $A({\bf k})$,
we first observe $A({\bf k})=A[\pi-1]$, and
the following Proposition.

\begin{prop}\label{PAr1}
For any ${\bf k}$-isogeny $r:A\rightarrow A$,
there is an isomorphism of $\zz[\pi]$-modules:
$A[r]\cong\prod_{l} T_l(A)/r T_l(A)$
where $l$ ranges over all prime numbers.
\end{prop}
\begin{proof}
The finite abelian group $A[r]$ has the decomposition
$A[r]\cong \prod_{l}A[r][l^\infty]$,
where each component is isomorphic to
$T_l(A)/r T_l(A)$ by the argument before
Proposition~\ref{P21}. All maps are $\zz[\pi]$-module homomorphisms.
\end{proof}

\subsection{Elementary supersingular abelian varieties}

It is well-known (see~\cite[Theorem~4.2]{Oort:74}) that
an abelian variety $A$ over ${\bf k}$
is supersingular if and only if either one of the following three conditions
holds:
(1) the eigenvalues of the Frobenius $\pi$ are supersingular $q$-numbers;
(2) the Newton polygon of $A$ is a straight line of slope $1/2$;
(3) $A$ is $\bar{\bf k}$-isogenous to a power of a supersingular
elliptic curve.

Note that $A[p]=0$ is the same as $T_p(A)=0$. We would like to clarify
the following facts without proof: A supersingular abelian variety $A$
over ${\bf k}$ has $A[p] = 0$ and the converse holds when $d=1$ or
$2$. However, the converse does not always hold when $d>2$.  In fact, an
abelian variety has $A[p]=0$ if and only if its Newton polygon has no
$0$-slope segment, which does not imply it being a straight line of slope
$1/2$ when $d>2$.

For the rest of this section we assume that $A$ is an elementary
supersingular abelian variety over ${\bf k}$ whose Frobenius relative
to ${\bf k}$ is $\pi$. The characteristic polynomial of $\pi$ is
$f=g^e$ for some monic irreducible polynomial $g$ over
$\qq$ and a positive integer $e$.
Since $\qq[\pi]=\qq(\pi)$ is a field, we fix an embedding of
$\qq(\pi)$ in $\cc$ and identify $\pi$ with its image, which is an
algebraic integer of the form $\zeta_m^{\nu}\sqrt{q}$ for some
primitive $m$-th root of unity $\zeta_m^\nu$ and the positive square
root $\sqrt{q}$.  We resume the notation from section~\ref{S3}, that
is, $K=\qq(\pi^2)=\qq(\zeta_{m/(2,m)})$, its ring of integers
$\O_K=\zz[\zeta_{m/(2,m)}]$, and $\O$ the ring of integers of
$\qq(\pi)$.

If given a supersingular $q$-number $\pi=\zeta_m^\nu\sqrt{q}$, we
describe the endomorphism algebra of $A$ over ${\bf k}$ in the sense
described on the proposition below. Let $\Q$ be the set of all
supersingular $q$-numbers $\zeta_m^\nu\sqrt{q}$ for some primitive
root of unity $\zeta_m^\nu$ such that either of the following two
conditions is satisfied: (1) $m=1$ or $2$; (2) $q$ is a square,
$(2,p)p\nmid m$ and $p$ is of odd order in the group
$(\zz/m_{(p)}\zz)^*$.

\begin{prop}\label{Pe12}
Suppose $A$ is simple supersingular over ${\bf k}$ with Frobenius
$\pi$.
\begin{enumerate}
\item[(1)] If $\pi\in\Q$ then $e=2$ and $\End_{\bf k}^0(A)$
is a quaternion algebra
over $\qq(\pi)$;
\item[(2)]
If $\pi\not\in\Q$ then $e=1$ and $\End_{\bf k}^0(A)$
is commutative and equal to $\qq(\pi)$.
\end{enumerate}
\end{prop}
\begin{proof} Let $v$ be any place of $\qq(\pi)$ (including
both finite and infinite primes).  Let $e_v$ denote the denominator of
the Hasse invariant, $\inv_v(\End_{\bf k}^0(A))$, of $\End_{\bf
k}^0(A)$ at $v$.  By~\cite[Th\'eor\`eme 1]{Tate1:68} we have
$$\inv_v(\End_{\bf k}^0(A))=\frac{\ord_v(\pi)[\qq(\pi)_v:
\qq_p]}{\ord_v(q)}=
\frac{[\qq(\pi)_v:\qq_p]}{2}=\frac{\gamma(v/p)\kappa(v/p)}{2} \bmod
1,$$ for all primes $v$ lying over $p$, so $e_v=1$ or $2$.  Now
$e_v=1$ for all complex $v$ and also for all finite primes $v$ not
lying over $p$, while $e_v=2$ for all real $v$. We have $e=\lcm_v(e_v)
= 2$ if either (1)' $v$ is real or (2)' $\gamma(v/\wp)\kappa(v/\wp)$
is odd; and $e=1$ otherwise.  It is obvious that (1)' is equivalent to
(1).  We show below that if $v$ is not real prime then (2)' is
equivalent to (2):

Suppose $q$ is not a square: we claim that $e_v=1$ for all finite
primes $v$ over $p$.  Now $[\qq(\pi):K]=1$ or $2$. The former implies
$2\mid \gamma(v/p)$. Consider the latter case, if $\sqrt{p}\in\qq(\pi)$,
then $2\mid \gamma(v/p)$ and so $e_v=1$; otherwise, we would have quadratic
extensions $\qq(\zeta_m,\sqrt{p}) \supset \qq(\pi)\supset K$.  But if $p$
was unramified in $\qq(\pi)/K$, then it would be unramified in
$\qq(\zeta_m,\sqrt{p})/\qq(\zeta_m)$, which is
absurd; so we must conclude that $p$ is totally ramified in
$\qq(\pi)/K$ and hence $2\mid \gamma(v/p)$ and so $e_v=1$.

Suppose $q$ is a square: so that
$\qq(\pi)=\qq(\zeta_m)$. Then for any finite prime $v$ over $p$, we have
that $\kappa(v/p)$ equals the order of $p$ in
$(\zz/m_{(p)}\zz)^*$; Let $\phi(\cdot)$
denote the Euler phi-function here, then $\gamma(v/p)=\phi(m_{(p)})$,
which is odd if and only if $(2,p)p\nmid m$.
This finishes our proof.
\end{proof}

\begin{remark}\label{Re12}
Suppose $A$ is simple supersingular over ${\bf k}$.  If $\pi\in \E$,
then $\pi\in \Q$ if and only if $d = 2$.  This follows from the above
proposition and the definitions of $\E$ and $\Q$.  The remark will be
used in the proof of Proposition~\ref{Psimple2} in the future.
\end{remark}

\begin{remark}\label{R1}
Let $A$ be a simple supersingular abelian variety with
odd dimension $d>2$, then $e=1$ and $\End_{\bf k}^0(A)$ must be commutative.
Indeed, recall~\cite[Th\'eor\`eme 1]{Tate1:68} that $2d=e[\qq(\pi):\qq]$
and so it suffices to show $2\mid [\qq(\pi):K][\qq(\zeta_{m/(2,m)}):\qq]$.
Either $[\qq(\pi):K]=1$ or $2$, in
the former case $[\qq(\zeta_{m/(2,m)}):\qq]=\phi(m/(2,m))>1$ and so is even.
\end{remark}

\subsection{Module structures}

Let $R$ be a subring in $\qq(\pi)$ with $\zz[\pi]\subseteq R\subseteq
\End_{\bf k}(A)\cap \qq(\pi)$. For any finite group $G$, we write
$\#G$ for its order.

\begin{lemma}\label{L52}
Let $M'\subseteq M''$ be modules over any ring $R$. Let $r\in R$ be
such that $R/rR$ is finite and $r$ acts faithfully on
$M', M''$.
\begin{enumerate}
\item[(1)] If $M'$ contains a free $R$-module of rank
$s$ as a submodule of finite index, then $\#M'/r M'=(\#(R/r R))^s.$
\item[(2)] If $M', M''$ contain a free $R$-module of rank $s$
as a submodule of finite index in $M',M''$, respectively,
then there are  homomorphisms
$\rho': M'/r M'\rightarrow M''/r M''$
and  {$\rho'': M''/r M''\rightarrow M'/r M'$}
with $$\#\kernel(\rho')=\#\cokernel(\rho')
=\#\kernel(\rho'')=\#\cokernel(\rho'')
\;\mid\;\# M''/M'.$$
\end{enumerate}
\end{lemma}
\begin{proof}
(1) Since $r$ acts faithfully on $M'$ and $R^s$, the injective map
$r: M'\rightarrow M'$ induces an injective map
$r: R^s\hookrightarrow R^s$.
On the other hand,  the given injection $R^s\hookrightarrow M'$
is of finite index, we thus have
$\#(M'/r M')\cdot\#(M'/R^s)=\#(R^s/r
R^s)\cdot\#(M'/R^s).$ Therefore, $\#M'/r M'=\#(R/r
R)^s$.

(2) Let $r$ act on the short exact sequence of $R$-modules
$0\rightarrow M'\rightarrow
M''\rightarrow M''/M'\rightarrow 0,$ and apply the Snake lemma.
We then get the desired map $\rho'$ with $\#\cokernel(\rho')$
dividing $\#M''/M'$.
By part (1), we have $\#M'/r M'=\#M''/r M''$ as they both equal
$\#(R/r R)^s$. Thus $\#\kernel(\rho')=\#\cokernel(\rho')$.
Any finite $R/r R$-module $N$ has an isomorphic dual
$\Hom_\zz(N,\qq/\zz)$, our assertion on $\rho''$ follows by
taking the dual of $\rho'$.
\end{proof}

\begin{prop}\label{PAr} 
Let $r$ be an isogeny in $R$. Then there
is an $R$-module homomorphism
$$\varphi_r: A[r]\longrightarrow \prod_{l\neq p}(R_l/r R_l)^e$$ which
is an isomorphism except when $\pi\in\E$ in which case
$\#\kernel(\varphi_r)$ and $\#\cokernel(\varphi_r)$ are equal and divide
$2^d_{(p)}.$
\end{prop}
\begin{proof}
By Proposition~\ref{PAr1} and the fact $A[p]=0$, we have
$A[r]\cong\prod_{l\neq p}T_l/rT_l$.  Recall that $T_l$ is a
torsion-free $R_l$-module of rank $e$, so we invoke
Proposition~\ref{Pindecomposable}.  If $\pi\not\in\E$ or $p=2$, then
$T_l/rT_l\stackrel{\sim}{\rightarrow} (R_l/r R_l)^e$ for each $l\neq
p$, and we obtain the desired isomophism $\varphi_r$. Now suppose
$\pi\in\E$.  Lemma~\ref{Lindex}~(2) implies
$\#\O_\wp/R_\wp\;\mid\;\#\O_\wp/\zz[\pi]_\wp =\#(\O_K)_\wp/\wp =
2^{\kappa(\wp/2)}.$ Clearly $\kappa(\wp/2)\varrho(\wp/2)\mid [K:\qq]$
and $[K:\qq]=[\qq(\pi):\qq]/2$ by Lemma~\ref{Lindex}~(1).  For each
$l$, we have a map $T_l/rT_l{\rightarrow}(R_l/rR_l)^e$ which is an
isomorphism if $l\neq 2$.  When $l=2$, Lemma~\ref{L52} indicates the
size of its kernel and cokernel are equal and divide $(\#
T_2/R_2^e)_{(p)}$. Taking product over all $l\neq p$ we obtain the
desired map $\varphi_r$ with
$\#\kernel(\varphi_r)=\#\cokernel(\varphi_r)$ and divides
$$(\#\O_2/R_2)_{(p)}^e
\;\mid\; 2^{\kappa(\wp/2) \varrho(\wp/2)e}_{(p)}
\;\mid\; 2^{e[K:\qq]}_{(p)}
\;\mid\; 2^{e[\qq(\pi):\qq]/2}_{(p)}$$
where the last number equals $2^d_{(p)}$.
\end{proof}

\noindent {\bf Proof of Theorem~1.2.} Let $S=\zz-p\zz$. By
Proposition~\ref{PAr}, there is an $R$-module homomorphism $\varphi_n:
A[n]\rightarrow ((\frac{1}{n}R)/R)^e$ for every $n\in S$. Let $W_n$ be
the set of such homomorphisms.  If {$m\mid n$}, then by passing to the
largest submodule annihilated by $m$ we see that any $R$-module homomorphism
$\varphi_n$ maps the submodule $A[m]$ of $A[n]$ to
$((\frac{1}{m})R/R)^e$, so there is a restriction map $W_n\rightarrow
W_m$. Since the projective limit of a system of non-empty finite sets
is non-empty, the projective limit of the sets $W_n$ is
non-empty. Therefore we can make a simultaneous choice of
$R$-module homomorphisms $\varphi_n$ that commute with the inclusions
$A[m]\subseteq A[n]$ and $((\frac{1}{m}R)/R)^e\subseteq
((\frac{1}{n}R)/R)^e$. Taking the injective limit over $n\in S$, we
get an $R$-module homomorphism $\injlim A[n]
\rightarrow\injlim((\tfrac{1}{n}R)/R)^e$, that is $\varphi:
A(\bar{\bf k})\rightarrow (R_{(p)}/R)^e.$ Since
$A(\bar{\bf k})$ and $(R_{(p)}/R)^e$ are both divisible as abelian groups,
the cokernel of $\varphi$ is also divisible, but it is finite and hence
trivial. So $\cokernel(\varphi)\cong\injlim\cokernel(\varphi_n)$ is
trivial and $\varphi$ is surjective.  In $A(\bar{\bf k})$ we have
$\kernel(\varphi)\cong\injlim\kernel(\varphi_n)$.  Thus $\varphi$ is
an isomorphism except when $\pi\in\E$, in which case
$\#\kernel(\varphi)$ divides $2^d_{(p)}$ since $\#\kernel(\varphi_n)$
divides $2^d_{(p)}$ for each $n$.  \qed

\begin{prop}\label{Psimple2}
Let $A$ be a simple supersingular abelian variety
over ${\bf k}$ with $f=g^e$. Let $R=\End_{\bf k}(A)\cap\qq(\pi)$.
If $p=2$ or $d\neq 2$, then
$A(\bar{\bf k})\cong_R (R_{(p)}/R)^e$.
If $p\neq 2$ and $d=2$, then there are
non-negative integers $a,b$ with $a+b=e$ and
$$A(\bar{\bf k})\cong_R (R_{(p)}/R)^a \times (\O_{(p)}/\O)^b.$$
\end{prop}
\begin{proof}
Let $\varphi: A(\bar{\bf k})\rightarrow (R_{(p)}/R)^e$ be defined as in
Theorem~1.2, which is an isomorphism unless $\pi\in\E$. Suppose
$\pi\in\E$. Then $\#\kernel(\varphi) = \#\cokernel(\varphi)$ is a
$2$-power.  Suppose $d\neq 2$.  Then $e=1$ by Remark~\ref{Re12} and so
$T_2$ is a torsion-free $R_2$-module of rank 1.  Recall that $K$ is a
cyclotomic field.  For any prime $\wp$ in $\O_K$ lying over $2$, write
$T_\wp$ for the $\wp$-adic completion of $T_2$ and $T_\wp:T_\wp=\{r\in
R_\wp\mid r T_\wp\subseteq T_\wp\}$.  Then $T_\wp:T_\wp = R_\wp$, so
$T_\wp$ is a fractional ideal of $R_\wp$.  Recall from
Lemma~\ref{Lindex} that $R_\wp$ is a Bass $(\O_K)_\wp$-order and thus
$T_\wp\cong_{R_\wp}R_\wp$ by~\cite[\S~2.6]{Buch-HWL:94}.  So
$T_2\cong_{R_2}R_2$ and this induces isomorphism
$T_2/2T_2\stackrel{\sim}{\rightarrow}R_2/2R_2$ by Lemma~\ref{L52}.
Thus $\varphi$ is an isomorphism.
Suppose $d=2$. Then $\pi\in\E$ implies $\pi=\pm\sqrt{(\frac{-1}{p})q}$
and $e=2$ by Remark~\ref{Re12}. In this case, $\wp=2$, so by
Proposition~\ref{Pindecomposable}, we have $A[n]\cong \prod_{l\neq p}T_l/nT_l
\cong (\frac{1}{n}R/R)^a\times
(\frac{1}{n}\O/O)^b$ for all $n\in S=\zz-p\zz$.
Take injective limit both sides over $n\in S$, we have
$$A(\bar{\bf k}) \cong \injlim((\tfrac{1}{n}R/R)^a\times
(\tfrac{1}{n}\O/\O)^b)\cong_R (R_{(p)}/R)^a\times
(\O_{(p)}/\O)^b.$$This finishes our proof.\end{proof}

\subsection{Group structures}

In this subsection we shall apply the results of the previous
subsection to our study of the group structure of $A({\bf k})$.

If $A$ is exceptional, $\qq(\pi)=\qq(\sqrt{(\tfrac{-1}{p})q})
=\qq(\sqrt{(\tfrac{-1}{p})p})$, so
$\O=\zz[\frac{1+\sqrt{(\frac{-1}{p})p}}{2}].$
By Lemma~\ref{Lindex}~(2) we notice
$\O_2/\zz[\pi]_2\cong \zz_2/2\zz_2\cong \zz/2\zz.$

\vspace{3mm}
\noindent {\bf Proof of Theorem~1.1.}
Apply Corollary~\ref{C58} to $M=T_l(A)$ and $R=\zz[\pi]$.
Now
$$A({\bf k})\cong_{\zz[\pi]}(\zz[\pi]/(\pi-1))^e
\cong_\zz (\zz/g(1)\zz)^e$$
unless $A$ is exceptional, in which case
the argument preceding the proof
implies that $(\pi-1)/2\in \O_2$ while $(\pi-1)/4\notin \O_2$. Since
$\#\O_2/(\pi-1)=\#\zz[\pi]_2/(\pi-1)=|g(1)|_2$, we have
$$\O_2/(\pi-1)\cong_{\zz_2}\zz_2/2\zz_2\times
\zz_2/\tfrac{g(1)}{2}\zz_2.$$ Hence there are non-negative integers
$a, b$ with $a+b=e$ such that
$$\begin{array}{lll} A({\bf
k})\!\!&\cong_{\zz[\pi]}&\!\!\!(\zz[\pi]_2/(\pi-1))^a\times(\O_2/(\pi-1))^b
\times \prod_{l\neq 2}(\zz[\pi]_l/(\pi-1))^e\\ \!\!&\cong_{\zz}
&\!\!\!  ((\zz_2/g(1)\zz_2)^a \times
(\zz_2/\tfrac{g(1)}{2}\zz_2\times\zz_2/2\zz_2)^b) \times\prod_{l\neq
2}(\zz_l/g(1)\zz_l)^{e}\\ \!\!&\cong_{\zz}
&\!\!\!(\zz/g(1)\zz)^{a}\times (\zz/\tfrac{g(1)}{2}\zz\times
\zz/2\zz)^{b}.
\end{array}$$
This proves our theorem.
\qed

\begin{prop}\label{Paprime}
Let the notation be as in Theorem~\ref{Telementary}. If $A$ is
exceptional, then for every pair of non-negative integers $a', b'$ with
$a'+b'=e$ there exists an abelian variety $A'$ isogenous over ${\bf k}$
to $A$ such that
$$A'({\bf k})\cong_\zz (\zz/g(1)\zz)^{a'}\times
(\zz/\tfrac{g(1)}{2}\zz \times \zz/2\zz)^{b'}.$$
\end{prop}
\begin{proof}
Let $A$ be exceptional. By Theorem~\ref{Telementary},
there are non-negative integers $a,b$ with
$a+b=e$ such that $T_2\cong_{\zz[\pi]_2}\zz[\pi]_2^{a}\times
\O_2^{b}$ and $$A({\bf k})\cong_\zz (\zz/g(1)\zz)^{a}\times
(\zz/\tfrac{g(1)}{2}\zz\times \zz/2\zz)^{b}.$$
If ${b'}=b$, then we are done.
If ${b'}<b$, let
$$M=\zz[\pi]_2^{a}\times\O_2^{b'}\times
(\tfrac{1}{2}\zz[\pi]_2)^{b-{b'}};$$
if $b'>b$, let
$$M=\zz[\pi]_2^{a'}\times\O_2^{b'},$$ in either case
$M\cong_{\zz[\pi]_2}\zz[\pi]_2^{a'}\times\O_2^{b'}$.  By the argument
preceding the proof of Theorem~1.1, we know that $\O_2\subseteq
\frac{1}{2}\zz[\pi]_2\subset \qq(\pi)_2$.  By Proposition~\ref{P21},
there exits an abelian variety $A'$ over ${\bf k}$ with $T_2(A')=M$
and a ${\bf k}$-isogeny $A\stackrel{r}{\rightarrow}A'$ with $A[r]\cong
T_2(A')/T_2(A)$ while $T_l(A')=T_l(A)$ for all $l\neq 2$. Thus
$$\begin{array}{lll}
A'({\bf k})&\cong&\prod_{l\neq p}T_l(A')/(\pi-1)T_l(A')\\
&\cong_{\zz[\pi]}&(\zz_2[\pi]/(\pi-1))^{a'}\times(\O_2/(\pi-1))^{b'}
\times \prod_{l\neq 2}(\zz[\pi]_l/(\pi-1))^e\\
&\cong_\zz& (\zz/g(1)\zz)^{a'}\times
(\zz/\tfrac{g(1)}{2}\zz \times \zz/2\zz)^{b'}.
\end{array}$$
This finishes the proof.
\end{proof}

\begin{corollary}\label{Csimple}
Suppose $A$ is a simple supersingular abelian variety over ${\bf k}$
of dimension $d>2$ with $f=g^e$, then $A({\bf k})\cong_\zz
(\zz/g(1)\zz)^e$ with $e=1$ or $2$.  If $d =1$, then $A$ is a
supersingular elliptic curve and $A({\bf k})\cong_\zz (\zz/g(1)\zz)^e$
or $A({\bf k})\cong_\zz \zz/\tfrac{q+1}{2}\zz \times \zz/2\zz$; that
latter case occurs only when $q$ is not a square and $p\equiv 3\bmod
4$.  If $d =2$, then $A$ is a simple supersingular abelian surface
and $A({\bf k})\cong_\zz (\zz/g(1)\zz)^e$ or $A({\bf k})\cong_\zz
\zz/\tfrac{q+1}{2}\zz \times \zz/2\zz$; that latter case occurs only
when $q$ is not a square and $p\equiv 1\bmod 4$.
\end{corollary}
\begin{proof}
If $A$ is simple over ${\bf k}$ of dimension $d>2$, then
$A$ is never exceptional, so $A({\bf k})\cong_\zz (\zz/g(1)\zz)^e$,
where $e=1$ or $2$ as we have seen in Proposition~\ref{Pe12}.

If $A$ is an elliptic curve, then
$A({\bf k})\cong_\zz (\zz/g(1)\zz)^e$ unless $A$ is exceptional in which case
$A({\bf k})\cong_\zz (\zz/g(1)\zz)^e$ or
$A({\bf k})\cong_\zz \zz/\tfrac{q+1}{2}\zz \times \zz/2\zz$.
Both cases may occur because of Proposition~\ref{Paprime}.
(This result can be found in~\cite[Chapter~4, (4.8)]{Schoof:87}.)

If $A$ is of dimension $2$, then
$A({\bf k})\cong_\zz(\zz/g(1)\zz)^2$
unless $A$ is exceptional in which case
$A({\bf k})\cong_\zz(\zz/g(1)\zz)^2$ or
$A({\bf k})\cong_\zz \zz/\frac{q-1}{2}\zz \times \zz/2\zz$.
\end{proof}

In particular, by Remark~\ref{R1}, if $A$ is simple supersingular of odd
dimension $d>2$, then $A({\bf k})\cong \zz/g(1)\zz$.


\begin{thebibliography}{99}

\bibitem{Buch-HWL:94}
{\sc J. Buchmann and H.W. Lenstra, Jr.:},
{Approximating rings of integers in number fields},
{\it J. Th\'eor. Nombres Bordeaux},{\bf 6}(1994), 221--260.

\bibitem{Curtis-Reiner:90}
{\sc C.W. Curtis and I. Reiner:},
{\it Methods of representation theory},
{J.~Wiley and Sons},
1990, volume 1.

\bibitem{Deuring:41}
{\sc M. Deuring:}
Die {T}ypen der {M}ultiplikatorenringe elliptischer
          {F}unktionenk\"{o}rper,
{\it Abh. Math. Sem. Univ. Hamburg.},
{\bf 14}(1941), 197--272.

\bibitem{Frohlich-Taylor:91}
{\sc A. Frohlich and M.J. Taylor:}
{\it Algebraic number theory},
Cambridge University Press, 1991.
Vol. 27,
Cambridge studies in advanced mathematics.

\bibitem{Greither:82}
{\sc C. Greither:}
On the two generator problem for the ideals of a one-dimensional ring,
{\it J. Pure Appl. Algebra}, {\bf 24} (1982), 265--276.

\bibitem{HWL:96}
{\sc H.W. Lenstra, Jr.:}
Complex multiplication structure of elliptic curves,
{\it J. Number Theory},
{\bf 56}(1996), 227--241.

\bibitem{Lang:86}
{\sc S. Lang:}
{\it Algebraic number theory},
Graduate Texts in Mathematics,
vol. 110, Springer-Verlag.

\bibitem{Levy:85}
{\sc L. Levy and R. Wiegand},
Dedekind-like behavior of rings with $2$-generators,
{\it J. Pure Appl. Algebra},
{\bf 37} (1985), 41--58.

\bibitem{Milne}
{\sc J. Milne:}
Ed. G. Cornell and J. Silverman,
Chapter V, Abelian varieties,
{\it Arithmetic geometry},
{Springer-Verlag}, 1987.


\bibitem{Mumford:74}
{\sc D. Mumford:}
{\it Abelian varieties},
Oxford University Press,
1974.

\bibitem{Oort:74}
{\sc F. Oort:}
Subvarieties of moduli spaces,
{\it Invent. Math.},
{\bf 24} (1974), 95--119.

\bibitem{Schoof:87}
{\sc R. Schoof:}
Nonsingular plane cubic curves over finite fields,
{\it J. Combin. Theory. Ser. A},
{\bf 46}(1987), 183--211.

\bibitem{Silverman:86}
{\sc J. Silverman:},
{\it The arithmetic of elliptic curves},
Springer-Verlag,
Graduate Texts in Mathematics,
{\bf 106}, 1986.

\bibitem{Tate1:68}
{\sc J. Tate:}
Classes d'isog\'enie
des vari\'et\'es ab\'eliennes sur un corps fini (d'apr\`es {T}.
{H}onda),
{\it S\'eminaire Bourbaki},
{\bf 21} (1968/69),
no. 352, 95--110.

\bibitem{Waterhouse:69}
{\sc W. Waterhouse:}
Abelian varieties over finite fields,
{\it Ann. Sci. \'Ecole Norm. Sup.},
{\bf 2}(1969), 521--560.

\bibitem{Waterhouse-Milne:71}
{\sc W. Waterhouse and J. Milne:}
Abelian varieties over finite fields,
{\it Proc. Sympos. Pure Math.},
{\bf 20}(1971),53--64.

\bibitem{Xing:96}
{\sc C-P. Xing:}
On supersingular abelian varieties of
                  dimension two over finite fields,
{\it Finite Fields Appl.} {\bf 2}(1996), 407--421.

\bibitem{paper2}
{\sc H.J. Zhu:}
Supersingular abelian varieties over finite fields,
{\it J. Number Theory}, {\bf 86} (2001), 61--77.
\end{thebibliography}
\end{document}